\documentclass[preprint,12pt]{elsarticle}
\usepackage{graphicx}
\usepackage{amsmath}
\usepackage{natbib}
\usepackage{amssymb}
\usepackage{hyperref}
\biboptions{}
\journal{Mathematics and Computers in Simulation}
\begin{document}
\begin{frontmatter}
\title{Meshfree method for fluctuating hydrodynamics}
\author[AP]{Anamika  Pandey}
\ead{pandey@mathematik.uni-kl.de}
\author[AP]{Axel Klar}
\ead{klar@mathematik.uni-kl.de}
\author[AP]{Sudarshan Tiwari}
\ead{tiwari@mathematik.uni-kl.de}
\address[AP]{Fachbereich Mathematik, TU Kaiserslautern, D-67653  Kaiserslautern, Germany}
\begin{abstract}
In the current study  a meshfree Lagrangian particle method for the Landau-Lifshitz Navier-Stokes (LLNS) equations is developed. The LLNS equations incorporate thermal fluctuation into macroscopic hydrodynamics by the addition of white noise fluxes whose magnitudes are set by a fluctuation-dissipation theorem.  The study focuses on capturing the correct variance and correlations computed at equilibrium flows, which are compared with available theoretical values. Moreover, a numerical test for the random walk of standing shock wave has been considered for capturing the shock location.
\end{abstract}
\end{frontmatter}
\section{Introduction}
\label{sec1}
Physical quantities which describe a macroscopic system in equilibrium seem to be very near to their mean value. Nevertheless, due to microscopic fluctuation, random deviation from this mean value, though small, do occur. Thermal fluctuation are a source of noise in many system. These fluctuation play a major role in phase transitions and chemical kinetics.  

Investigating thermal fluctuation in the motion of fluids becomes essential at micro and nano scale, because of the various applications of micro and nano scale flow, ranging from micro-engineering to molecular biology. Micro-machines have a major impact on many disciplines (e.g. biology, medicine, optics, aerospace, and mechanical and electrical engineering) \cite{Karniadakis, C.M. Ho}. 

The study of fluctuation at micro and nanoscale is particularly interesting when the fluid is under extreme conditions or near a hydrodynamic instability, e.g. the breakup of droplet in nanojet, fluid mixing in the Rayleigh-Taylor instability \cite{moseler, eggers}. 

The presence of thermal fluctuation becomes significant for larger Kundsen number $(Kn \ge 0.01)$\footnote{$Kn = \frac{\lambda}{L}$, $\lambda$ denotes mean free path and $L$ represents the characteristic length.} The LLNS equations  try to capture these thermal fluctuation as accurately  as possible which is not possible in the case of Navier-Stokes equations. 

 To describe the general theory of fluctuation in fluid dynamics is equivalent to setting up the "equation of motion" for fluctuating quantities. Landau and Lifshitz introduced the appropriate additional terms in the general equation of fluid dynamics and gave an extended form of the Navier-Stokes equations. The Landau-Lifshitz Navier-Stokes equations are written as
 \begin{equation}
 \label{eq1}
 \textbf{U}_{t} + \nabla . \textbf{F} = \nabla . \textbf{D} + \nabla . \textbf{S},
 \end{equation} 
 where, U stands for the vector of conserved quantities, density of mass, momentum and energy 
\begin{equation}  
\label{eq2}
 \textbf{U} = \left( \begin{array}{c} \rho \\ \textbf{J} \\ E \end{array}\right),
 \end{equation}
  $\textbf{F}$ denotes the hyperbolic flux and $\textbf{D}$ denotes the diffusive flux of fluid dynamic equations. $\textbf{F}$ and $\textbf{D}$ are given by
 \begin{equation}
 \label{eq3}
 \textbf{F} = \left( \begin{array}{c} \rho \textbf{v} \\ \rho \textbf{v}\textbf{v} + P \textbf{I}\\ \textbf{v}E + P \textbf{v}\end{array} \right),
 \end{equation}
 \begin{equation}
 \label{eq4}
 \textbf{D} = \left( \begin{array}{c} 0 \\ \tau \\ \tau .\textbf{v} - \textbf{q} \end{array}\right),
 \end{equation}
 where $\textbf{v} $ is the fluid velocity, $P$ is the pressure and $T$ denotes the temperature. $\tau = \eta\left( \nabla \textbf{v} + \nabla \textbf{v}^{T} - \dfrac{2}{3}\textbf{I}\nabla \cdot \textbf{v}\right) $ is the stress tensor. $\textbf{q} =  - \kappa \nabla T $ denotes the heat flux. Here $\eta$ and $\kappa$ are the coefficients of viscosity and thermal conductivity, respectively. For the given expression of $\tau$ we have assumed the bulk viscosity to be zero.
 
 The  expression for $\tau$ and $\textbf{q}$ relate these quantities to the velocity and temperature gradients respectively. But, in the presence of fluctuation there are also spontaneous local stresses and heat fluxes in the fluid, which are not related to velocity and temperature gradient. For these spontaneous local stresses tensor and heat fluxes, the LLNS equations introduce additional quantities in the fluid dynamic equations called stochastic flux, i.e.
 \begin{equation}
 \label{eq5}
 \textbf{S} = \left( \begin{array}{c} 0 \\ \textit{S} \\ \textit{H} + \textbf{v} \cdot \textit{S} \end{array}\right),
 \end{equation}
where the stochastic stress tensor (sst) \textit{S} and stochastic heat flux (shf) \textit{H} have zero mean and their covariances are given by
\begin{equation}
\label{eq6}
Cov(\textit{S}_{ij}(\textbf{r},t),\textit{S}_{kl}(\textbf{r}^{'},t^{'})) = 2k_{B} \eta T \left(\delta_{ik}^{K}\delta_{jl}^{K} + \delta_{il}^{K} \delta_{jk}^{K} - \frac{2}{3}\delta_{ij}^{K}\delta_{kl}^{K}\right)\delta(\textbf{r}-\textbf{r}^{'})\delta (t-t^{'}),
\end{equation}
\begin{equation}
\label{eq7}
Cov(\textit{H}_{i}(\textbf{r},t),\textit{H}_{j}(\textbf{r}^{'},t^{'})) = 2k_{B}\kappa T^{2}\delta_{ij}^{K}\delta(\textbf{r}-\textbf{r}^{'})\delta(t-t^{'}),
\end{equation}
\begin{equation}
\label{eq8}
Cov(\textit{S}_{ij}(\textbf{r},t),\textit{H}_{k}(\textbf{r}^{'}, t^{'})) = 0,
\end{equation}
where, $k_{B}$ is the Boltzmann's constant. 

These stochastic properties for $\textit{S}$ and $\textit{H}$ have been derived by a variety of approaches. Originally, these properties have  been derived for equilibrium fluctuation \cite{Landau1, bixon, fox, kelly} and later the validity of the LLNS equations for non-equilibrium systems has been shown \cite{espa}. 

In this work a meshfree numerical scheme is developed for solving the LLNS equations. For simplicity, we will deal with one-dimensional system.  The Lagrangian form of  the LLNS equations for $1D$ system in terms of primitive variables can be written as 
\begin{equation}
\label{eq9}
\dfrac{D\rho}{Dt} = -\rho \dfrac{\partial u}{\partial x}
\end{equation}
\begin{equation}
\label{eq10}
\rho \dfrac{Du}{Dt} =  - \dfrac{\partial P}{\partial x} + \dfrac{\partial}{\partial x}\left( \dfrac{4}{3} \eta \partial_{x}u \right) + \dfrac{\partial \textit{s}}{\partial x}
\end{equation}
\begin{equation}
\label{eq11}
c_{v}\rho\dfrac{DT}{Dt} = - P \dfrac{\partial u}{\partial x} + \dfrac{4}{3}\eta\left( \dfrac{\partial u}{\partial x}\right)^{2} + \dfrac{\partial}{\partial x}\left( \kappa \dfrac{\partial T}{\partial x}\right)  + \textit{s}\dfrac{\partial u}{\partial x} + \dfrac{\partial \textit{h}}{\partial x}.
\end{equation}
$u$ is the fluid velocity in x-direction and $T$ is the temperature. $\textit{s}$ and $\textit{h}$ represent sst and shf in 1D respectively. Momentum $\textbf{J} = \rho u$ and energy density $E = \textit{c}_{v}\rho T + \frac{1}{2}\rho u^{2}$ is expressed in terms of $\rho, u, T$. By $D/Dt$ we denote the Lagrangian derivative. We will take the above system with equation of state $P = \rho R T$, where $R$ is the gas constant. We will demonstrate our result for a mono-atomic,  hard sphere gas for which $R = k_{B}/m $ and $ \textit{c}_{v} = \dfrac{R}{\gamma - 1}$ where $m$ is molecular mass and  $\gamma  (=\frac{5}{2})$  is the ratio of specific heat.

Now for the covariances of the stochastic fluxes in a $1D$ system one obtains  from equations (\ref{eq6}), (\ref{eq7}) and (\ref{eq8})
\begin{align}
\label{eq12}
Cov(\textit{s}(x,t),\textit{s}(x^{'},t^{'})) &= \dfrac{1}{\sigma^{2}}\int dy \int dy^{'}\int dz \int dz^{'}Cov(\textit{S}_{xx}(\textbf{r},t), \textit{S}_{xx}(\textit{r}^{'},t^{'})) \nonumber \\ 
& = \dfrac{8k_{B}\eta T}{3\sigma}\delta(x-x^{'})\delta(t-t^{'})
\end{align}
and similarly,
\begin{align}
\label{eq13}
Cov(\textit{h}(x,t),\textit{h}(x^{'},t^{'})) &= \dfrac{1}{\sigma^{2}}\int dy \int dy^{'}\int dz \int dz^{'}Cov(\textit{H}_{x}(\textbf{r},t), \textit{H}_{x}(\textit{r}^{'},t^{'})) \nonumber \\ 
& = \dfrac{2k_{B}\kappa T^{2}}{\sigma}\delta(x-x^{'})\delta(t-t^{'}).
\end{align}
Here $\sigma$ represents the surface area of the system in the  $yz$ - plane.

A number of numerical schemes have been developed for stochastic hydrodynamic equations. A stochastic lattice-Boltzmann model has been developed for simulating solid-fluid suspensions by Ladd \cite{ladd}.  For modelling the Brownian motion of particles a similar approach has been used by Sharma and Patankar  in \cite{sharma}, where they coupled the fluctuating hydrodynamics equations with  equations of motion for the particle. Moseler amd Landman \cite{moseler} have used LLNS stochastic stress tensor for   lubrication equations and obtain good comparison with molecular dynamics simulation for the breakup of nanojets. 

Serrano and Espa\~nol \cite{serrano} have developed a thermodynamically consistent mesoscopic fluid particle model by casting their model into the so called GENERIC structure which allows to introduce thermal fluctuation. They describe a finite volume Lagrangian discretization of the continuum equations of hydrodynamics using Voronoi tessellation.  A similar mesoscopic, Voronoi-based algorithm using the dissipative particle dynamics method has been used by Fabritiis et al.\cite{fabr}.

Garcia et al.\cite{Garcia1} have developed a simple finite difference scheme for the linearized LLNS equations. This scheme has been designed for specific problems. In the context of adaptive mesh and algorithm refinement hybrid schemes  that couple continuum and particle algorithms have been developed, see, for example,  \cite{garcia5},  \cite{garcia6} and  \cite{garcia7}. Further, an algorithm refinement and coupling   molecular dynamics simulations to numerics of stochastic hydrodynamic equations have been presented by Fabritiis \cite{fabr1}.

In a paper by  Garcia et al. \cite{Garcia2} CFD based schemes for stochastic PDEs have been developed, where numerical scheme for the full LLNS equations have been developed. Spatial and time correlation at equilibrium have  been compared with theoretical values and DSMC simulations. Moreover, the effect of fluctuations on the shock drift has been shown and results are  compared with DSMC simulation.  The method is based on a third order, TVD Runge-Kutta temporal integrator (RK3) combined with a centered discretization of hyperbolic and diffusive fluxes. This scheme  also incorporates a specific interpolation for the required accuracy in variance. 

In this work we will present a grid free method for the LLNS equations. We will consider  a particle method  based on a least squares approach, see \cite{tiwari1}. We concentrate on  capturing the  correct variance in equilibrium flow and compare the result with theoretical values. We will also show the fluctuation effect in the shock location for a standing shock wave, mentioned above. The concluding section will discuss future work. Since, the developed method is a grid free method and the distribution of particles (moving grid) can be quite arbitrary, the method is suitable for complicated geometry and multiphase flows, see for example,  \cite{tiwari3} for the classical case without stochastic fluxes.
\section{Numerical Method}
To extend the ideas of meshfree methods for the  Navier-Stokes equations to the LLNS equations we will consider the  $1D$ LLNS equations in a Lagrangian framework.  To develop a meshfree framework for stochastic partial differential equations (SPDE) has a number of advantage for example, for the  study of the dynamics of small particles at fluid interfaces or for the development of  hybrid methods for coupling Monte Carlo methods for the  Boltzmann equation and the fluctuating hydrodynamics equations. We mention that earlier studies have dealt with  the coupling of DSMC for the Boltzmann equation and finite volume methods for the LLNS \cite{garcia8}.  

In the Lagrangian framework  an additional equation for the particle position, i.e.
\begin{equation}
\label{eq14}
\dfrac{Dx}{Dt} = u
\end{equation} is solved together  with (\ref{eq9} - \ref{eq11}). 
Here $u$ is the fluid velocity and $x$ denote the position of particle in $1D$. To approximate spatial derivatives at every grid point is equivalent to approximate the spatial derivative at every particle positions. For solving the Lagrangian LLNS system given by (\ref{eq9}-\ref{eq11}) together with (\ref{eq14}), we fill the domain by particles,\footnote{initially these particle are just a regular grid  with equal spacing}. These particles move with fluid velocity. Then  the  spatial derivative in equation (\ref{eq9}-\ref{eq11}) are approximated at each particle position from its neighbouring particles. This reduces the given system of stochastic partial differential equations (SPDEs) to a system of stochastic ordinary differential equations (SODEs) with respect to time. 

We will use the MacCormack scheme \cite{Garcia2} for the system of SODEs. The discretized form of the  LLNS equations is
\paragraph*{Step 1}
\begin{align}
\label{eq15}
x_{i}^{*} &= x_{i}^{m} + \triangle t  u_{i}^{m}\\
\rho_{i}^{*} &= \rho_{i}^{m} - \triangle t  \rho_{i}^{m}\left( \dfrac{\partial u}{\partial x}\right) ^{m}_{i} \\
u_{i}^{*} &= u_{i}^{m} + \dfrac{\triangle t}{\rho_{i}^{m}}\left\{ - \left( \dfrac{\partial P}{\partial x} \right) ^{m}_{i} + \dfrac{4}{3}\eta_{i}^{m}\left(  \dfrac{\partial^{2}u}{\partial x^{2}} \right) _{i}^{m}
+ \dfrac{4}{3}\left( \dfrac{\partial \eta}{\partial x}\right) _{i}^{m} \left( \dfrac{\partial u}{\partial x}\right) ^{m}_{i} + \left( \dfrac{\partial \textit{s}}{\partial x}\right) ^{m}_{i}\right\} \\
T_{i}^{*} &= T_{i}^{m} + \dfrac{\triangle t}{c_{v}\rho_{i}^{m}} \left\{ -P_{i}^{m} \left( \dfrac{\partial u} {\partial x} \right) ^{m}_{i} + \dfrac{4}{3}\eta_{i}^{m} \left( \left( \dfrac{\partial u}{\partial x}\right) ^{2} \right) ^{m}_{i} + \kappa_{i}^{m} \left( \dfrac{\partial^{2}T}{\partial x^{2}} \right) ^{m}_{i}\right. \nonumber \\ &+ \left. \left( \dfrac{\partial \kappa}{\partial x}\right) ^{m}_{i}\left( \dfrac{\partial T}{\partial x} \right) ^{m}_{i}
+ \textit{s}^{m}_{i}\left( \dfrac{\partial u}{\partial x}\right) _{i}^{m} 
+ \left( \dfrac{\partial \textit{h}}{\partial x}\right) ^{m}_{i}\right\}
\end{align}
\paragraph*{Step 2}
\begin{align}
x_{i}^{**} &= x_{i}^{*} + \triangle t  u_{i}^{*}\\
\rho_{i}^{**} &= \rho_{i}^{*} - \triangle t \rho_{i}^{*}\left( \dfrac{\partial u}{\partial x}\right) ^{*}_{i}\\
u_{i}^{**} &= u_{i}^{*} + \dfrac{\triangle t}{\rho_{i}^{*}}\left\{ - \left( \dfrac{\partial P}{\partial x} \right) ^{*}_{i} + \dfrac{4}{3}\eta_{i}^{*}\left(  \dfrac{\partial^{2}u}{\partial x^{2}} \right) _{i}^{*}
+ \dfrac{4}{3}\left( \dfrac{\partial \eta}{\partial x}\right) _{i}^{*} \left( \dfrac{\partial u}{\partial x}\right) ^{*}_{i} + \left( \dfrac{\partial \textit{s}}{\partial x}\right) ^{*}_{i}\right\}\\
T_{i}^{**} &= T_{i}^{*} + \dfrac{\triangle t}{c_{v}\rho_{i}^{*}} \left\{- P_{i}^{*} \left( \dfrac{\partial u} {\partial x} \right) ^{*}_{i} + \dfrac{4}{3}\eta_{i}^{*} \left( \left( \dfrac{\partial u}{\partial x}\right) ^{2} \right) ^{*}_{i} + \kappa_{i}^{*} \left( \dfrac{\partial^{2}T}{\partial x^{2}} \right) ^{*}_{i}\right. \nonumber \\ &+ \left. \left( \dfrac{\partial \kappa}{\partial x}\right) ^{*}_{i}\left( \dfrac{\partial T}{\partial x} \right) ^{*}_{i}
+ \textit{s}^{*}_{i}\left( \dfrac{\partial u}{\partial x}\right) _{i}^{*} 
+ \left( \dfrac{\partial \textit{h}}{\partial x}\right) ^{*}_{i}\right\}
\end{align}
\paragraph*{Final Step}
\begin{equation}
\label{eq23}
x^{m+1}_{i} = \dfrac{1}{2}\left( x^{m}_{i} + x_{i}^{**}\right) 
\end{equation}
\begin{equation}
\label{eq24}
\rho^{m+1}_{i} = \dfrac{1}{2}\left( \rho^{m}_{i} + \rho_{i}^{**}\right) 
\end{equation}
\begin{equation}
\label{eq25}
u^{m+1}_{i} = \dfrac{1}{2}\left( u^{m}_{i} + u_{i}^{**}\right) 
\end{equation}
\begin{equation}
\label{eq26}
T^{m+1}_{i} = \dfrac{1}{2}\left( T^{m}_{i} + T_{i}^{**}\right) 
\end{equation}

For each of the above steps  $P$, $\eta$, $\kappa$ will be computed by
\begin{equation}
\label{eq27}
P = \rho R T,
\end{equation}
\begin{equation}
\label{eq28}
\eta = \dfrac{5}{16\textit{d}^{2}}\sqrt{\dfrac{Mk_{B}}{\pi}T}.
\end{equation}
\begin{equation}
\label{eq29}
\kappa = \dfrac{15k_{B}\eta}{4M},
\end{equation}
Here, $\textit{d}$ denotes the molecular diameter, $M$ is molecular mass, $m = 0,1,2,\hdots$ represents the time step and $ N$ is the number of particles in the domain.

The approximation for sst and shf for each particle at any instant is computed as 
\begin{equation}
\label{eq30}
\textit{s}^{m}_{i} = \sqrt{\dfrac{8k_{B}}{3\triangle t V_{c}}\left( \eta_{i}^{m} T_{i}^{m}\right) }\hspace{2 mm} \Re^{m}_{i}
\end{equation}
\begin{equation}
\label{eq31}
\textit{h}^{m}_{i} = \sqrt{\dfrac{2k_{B}}{\triangle t V_{c}}\left( \kappa_{i}^{m} \left(T_{i}^{m}\right)^{2}\right) } \hspace*{2 mm}\Re^{m}_{i}
\end{equation}
where $V_{c}$ denotes the volume between two particle  and $\Re$ are independent, identically distributed (iid), Gaussian random variables with zero mean and unit variance.  

The stochastic fluxes require some extra care in multi-step scheme, variance in the stochastic flux $\textbf{\textit{s}} = (\textit{s}, \textit{h})$ is given by, 
\begin{align}
\label{eq32}
Var\left(\left(\textbf{\textit{s}}^{m+1}\right)^{2}\right) &= Var\left( \left( \dfrac{1}{2} \textbf{\textit{s}}^{m} + \dfrac{1}{2} \textbf{\textit{s}}^{*}\right)^{2}\right)  
 \nonumber \\ &= \left( \dfrac{1}{2} \right) ^{2} Var\left( \textbf{\textit{s}}^{m}\right)  + \left( \dfrac{1}{2}\right) ^{2}Var\left( \textbf{\textit{s}}^{*}\right)  \nonumber  \\
&= \dfrac{1}{2}Var\left( \textbf{\textit{s}}^{m}\right)  
\end{align}
on neglecting the multiplicity of noise i.e. $Var\left(\textbf{\textit{s}}^{m}\right) = Var\left(\textbf{\textit{s}}^{*}\right)$. 

Because of the temporal averaging, the  variance in the flux is reduced to half of its original magnitude. So, to include this observation the correct stochastic flux for a two step scheme will be$ \textbf{\textit{ \~s}} = \sqrt{2}\textbf{\textit{s}}$ instead of \textbf{\textit{s}}\cite{Garcia2}.

Now we have to solve  equations (\ref{eq15} - \ref{eq29}). The remaining task is to approximate the spatial derivatives on the right hand side of the prescribed equations. 
\subsection{Meshfree approximation of spatial derivatives}
We will describe the least square approximation of spatial derivatives in $1D$. As mentioned earlier, in this method grid points are particle positions. Therefore, we have to approximate the derivatives at every particle position. Let $f(t,x)$ be a scalar function at $x$ and $f_{i}(t)$ its value at $x_{i} \in [0, L]$ for $i = 1,2,3, \hdots, N$ for any instant $t$.  Spatial derivatives of $f(x)$ at $x$ will be approximated  on a set of neighbouring points. For limiting number of neighbouring points of $x$, one considers a weight function $w = w(\|x_{i} -x\| ; h)$ with small compact support, where $h$ determines the size of the support. We will consider a Gaussian weight function in the following form 
\begin{equation}
\label{eq33}
w\left(x_{i} - x ; h\right) =  \left\{ \begin{array}{ll} \exp\left( -\alpha\dfrac{\Vert x_{i} - x \Vert^{2}}{h^{2}}\right),  & \textrm{if $\dfrac{\Vert x_{i} - x \Vert}{h} \leq 1$} \\ 0, & \textrm{else} 
\end{array}
\right.
\end{equation}
with $\alpha$ a positive constant, chosen as  $\alpha=6.25$. $h$ defines the neighbourhood radius for $x$. Let $P(x,h) = \left\lbrace x_{i} : i = 1,2, \hdots , n\right\rbrace $ be the set of $n$ neighbouring points of $x$ in an interval of radius $h$. We have chosen  $h = 3dx,$ where $dx$ is the initial spacing of particles. 

Suppose we want to approximate the derivatives of a function $f(t,x)$ from its $n$ neighbouring points sorted with respect to its distance from $x$. Consider Taylor's expansion of $f(t,x_{i})$ around $x$
\begin{align}
\label{eq34}
f\left( t, x_{i}\right) = f\left( t, x\right) &+ \left( f\left( t, x\right) \right) _{x}\left(x_{i} - x \right)  \nonumber \\
&+\dfrac{1}{2}\left( f\left( t, x\right) \right)_{xx}\left( x_{i} - x\right)^{2} + e_{i}
\end{align}
where $e_{i}$ is the error in Taylor's expansion at the point $x_{i}$. The unknowns $f_{x}, f_{xx}$ are  computed by minimizing the error $e_{i}  \textrm{ for $ i = 1,2,3,\hdots,n$}$. The above system can be written as 
\begin{equation}
\label{eq35}
\vec{e} = M\vec{a}  - \vec{b}
\end{equation}
where,
\begin{equation}
\label{eq36}
M = \left( \begin{array}{cc} x_{1} - x & \dfrac{1}{2}\left( x_{1}-x\right)^{2}\\ [0.5 cm] x_{2} - x & \dfrac{1}{2}\left( x_{2}-x\right)^{2}\\[0.5 cm] \vdots & \vdots \\ [0.5 cm] x_{n} - x & \dfrac{1}{2}\left( x_{n}-x\right)^{2}.
\end{array}\right) ,
\end{equation}
$a = \left[ f_{x},f_{xx}\right] ^{T}\textrm{, $b =\left[f_{1}-f, f_{2}-f, \hdots, f_{n}-f\right]^{T}$ and $e = \left[e_{1},e_{2},e_{3}, \hdots, e_{n}\right]^{T}$.}$ 

For $n > 2$ this system will be over-determined for two unknowns $f_{x}$ and $f_{xx}$.

The unknowns $\vec{a}$ are obtained from a weighted least square method by minimizing the quadratic form
\begin{equation}
\label{eq37}
\mathit{J} =  \sum_{i = 1}^{n}w_{i}e_{i}^{2} = \left( M\vec{a} - \vec{b}\right)^{T} W \left(M\vec{a} - \vec{b}\right) 
\end{equation} 
where
\begin{displaymath}
W = \left(\begin{array}{cccc} w_{1} & 0 & \hdots & 0 \\ 0 & w_{2} & \hdots & 0 \\ \vdots & \vdots & \ddots & \vdots \\ 0 & 0 & \hdots & w_{n}
\end{array}\right)
\end{displaymath}
The minimization of $\mathit{J}$ gives
\begin{equation}
\label{eq38}
\vec{a} = \left(M^{T} W M\right)^{-1}\left(M^{T}W\right)\vec{b}
\end{equation}
and, finally, the required derivatives of  $f(t,x)$  as a linear combination of the discrete values $f_{i}$.
\section{Numerical Results}
\subsection{Equilibrium}
This section gives the results of the above described  method for an equilibrium scenario. The physical domain has been chosen such that the fluctuation in the system become significant.  The  parameters for the numerical simulation are given in Table \ref{table1}.  
\begin{table}[h]
\begin{tabular}{ | p {3 cm} | l |l| p {3 cm} | l | }
\hline
Molecular diameter (Argon) & $3.66 \times 10^{-8}$ && Molecular mass (Argon) & $6.63 \times 10^{-23}$ \\ [0.7 cm]
Reference mass density & $1.78 \times 10^{-3}$ && Reference temperature & $273$ \\ [0.7 cm]
Sound speed & $30781$ && Reference velocity & $0.5 \times 30781$ \\ [0.7 cm]
System length (L) & $1.25 \times 10^{-4}$ && Reference mean free path & $6.26 \times 10^{-6}$ \\ [0.7 cm]
System volume & $1.96 \times 10^{-16}$ && Time step & $1.0 \times 10^{-13}$ \\[0.7 cm] 
\hline
\end{tabular}
\caption{System parameter in CGS units for simulation of a dilute gas}
\label{table1}
\end{table}

The initial spacing of particle will be $dx = L /N$, where $N$ is the total number of initial particles. Here,  $N=40$  particles are  considered. This number is not  fixed during the simulation, particles are added or removed during the simulation. 

The stability condition is found to be consistent with that suggested in \cite{Garcia2}.
\begin{eqnarray}
\label{eq39}
\left(\mid u\mid + c_{s}\right)\dfrac{\triangle t}{\triangle x} \leq 1\\
max\left(\dfrac{4}{3}\dfrac{\bar{\eta}}{\bar{\rho}}, \dfrac{\bar{\kappa}}{\bar{\rho }c_{v}}\right) \dfrac{\triangle t}{\triangle x^{2}} \leq \dfrac{1}{2}
\end{eqnarray}
where $c_{s}$ is the sound speed, the bar indicates the reference value of quantities around which the system fluctuates. For the given reference state in Table \ref{table1} and the given  initial spacing of particles, the time step has been chosen as $\triangle t = 10^{-13} s$.
\subsubsection{Variance at Equilibrium}
The first benchmark  is to capture the  correct variance for fluctuations of  the  system at equilibrium. We consider a periodic domain with zero net flow and constant non-zero net flow. We take constant average density and temperature in both cases as given in Table \ref{table1}. The variance is computed from $10^{7}$ samples. We will calculate statistics in a global sense as given below.
\begin{align}
\label{eq40}
mean\left( \rho \right) &= \mathbf{E}(\rho) =  \dfrac{1}{\sum_{n = 1}^{N_{s}}M(n)} \left(\sum_{n = 1}^{N_{s}}\sum_{i = 1}^{M(n)}\rho_{i}^{n}\right), \\
Var\left( \rho \right) &= \mathbf{E}(\rho^{2}) - \left(\mathbf{E}(\rho)\right)^{2} \nonumber \\ &= \dfrac{1}{\sum_{n = 1}^{N_{s}}M(n)} \left(\sum_{n = 1}^{N_{s}}\sum_{i = 1}^{M(n)}\left(\rho_{i}^{n}\right)^{2}\right) \nonumber \\
&- \left(\dfrac{1}{\sum_{n = 1}^{N_{s}}M(n)} \left(\sum_{n = 1}^{N_{s}}\sum_{i = 1}^{M(n)}\rho_{i}^{n}\right)\right)^{2},
\end{align}
Here, $N_{s}$ is the total number of samples and $M(n)$ is the number of  particle at time $(nSkip + n)\triangle t$ ($nSkip$ is the number of initial time steps for stabilizing the system).
In the same way statistics for momentum and energy can be computed. 

Table \ref{table2} and \ref{table3} compare the theoretical variances which have been computed in (\cite{Landau2}, \cite{Garcia1}) with measured variances from the  meshfree stochastic scheme. 
\begin{table}[h]
\begin{tabular}{ | l | l | p {3 cm} | l |}
\hline 
\textbf{Variance of} & \textbf{Exact value} & \textbf{MacCormack with Meshfree} & \textbf{percentage error}  \\ [0.1 cm]
Density $\left(\rho\right)$ & $2.35 \times 10^{-8}$ & $2.11\times 10^{-8}$ & $ -10.2 \% $  \\ [0.1 cm]
Momentum $(\mathbf{J})$ & $13.34 $ & $13.33$ & $-0.07 \%$ \\ [0.1 cm]
Energy $(\mathbf{E})$ & $2.84 \times 10^{10}$ & $2.68\times 10^{10}$ & $-5.6 \% $\\ [0.1 cm]
\hline
\end{tabular}
\caption{Variance in conserved quantities at equilibrium for zero net flow.}
\label{table2}
\end{table}
\begin{table}[h]
\begin{tabular}{ | l | l | p {3 cm} | l |}
\hline 
\textbf{Variance of} & \textbf{Exact value} & \textbf{MacCormack with Meshfree} & \textbf{percentage error}  \\ [0.1 cm]
Density $\left(\rho\right)$ & $2.35 \times 10^{-8}$ & $2.12\times 10^{-8}$ & $ -9.7 \% $  \\ [0.1 cm]
Momentum $(\mathbf{J})$ & $18.91 $ & $19.01$ & $ +0.05 \%$ \\ [0.1 cm]
Energy $(\mathbf{E})$ & $3.67 \times 10^{10}$ & $3.85\times 10^{10}$ & $+4.6 \% $\\ [0.1 cm]
\hline
\end{tabular}
\caption{Variance in conserved quantities at equilibrium for constant non-zero net flow.}
\label{table3}
\end{table}
\subsubsection{Time covariance at equilibrium}
We will measure the time covariance of density fluctuation. The analytical formula of the time covariance for the density can be written as \cite{Garcia2}
\begin{align}
\label{eq43}
Cov\left( \rho(\omega , t) \rho(\omega , t + \tau)\right) = &
\lbrace\left( 1 - \dfrac{1}{\gamma}\right) exp\left\lbrace - \omega^{2}D_{T} \tau\right\rbrace + \dfrac{1}{\gamma}exp\left\lbrace - \omega^{2}\Gamma \tau\right\rbrace\cos \left(c_{s}\omega\tau\right) \nonumber \\
&+\dfrac{3\Gamma - D_{v}}{\gamma^{2}c_{s}}\omega exp \left\lbrace -\omega^{2} \Gamma\tau \right\rbrace \sin (c_{s}\omega\tau)\rbrace* Var\left( \rho(\omega , t)\right)
\end{align}
where $\omega = 2 \pi n/L $ is the wave number, $\gamma$ is the ratio of specific heat, $D_{T} = \kappa/\bar{\rho}c_{v}$ is the thermal diffusivity, $D_{v} = \dfrac{4}{3} \eta/\bar{\rho}$ is longitudinal kinematic viscosity, $c_{s}$ is the sound speed and $\Gamma = \dfrac{1}{2} \left[D_{v} + (\gamma -1)D_{T}\right]$ is the sound attenuation coefficient. 

For numerical simulation the time covariance is estimated from the mean of $N_{s}$ samples,
\begin{equation}
\label{eq44}
Cov\left(\rho(\omega, t) \rho(\omega, t + \tau)\right)_{N_{s}} =\dfrac{1}{N_{s}} \sum_{n = 1}^{N_{s}}R(t)R(t + \tau)
\end{equation}
\begin{equation}
Var\left(\rho\left(\omega, t\right)\right) = Var\left(R\left( t \right)\right)
\end{equation}
 where,
 \begin{equation}
 \label{eq44a}
 R(t) = \dfrac{1}{M(n)}\sum_{i = 1}^{M(n)}\rho_{i}\sin(2 \pi n x_{i}/L)
 \end{equation}
 We are considering the lowest wave number, i.e. $ n =1$.
\begin{figure}[h]
\centering
\includegraphics[width=8cm,height=8cm]{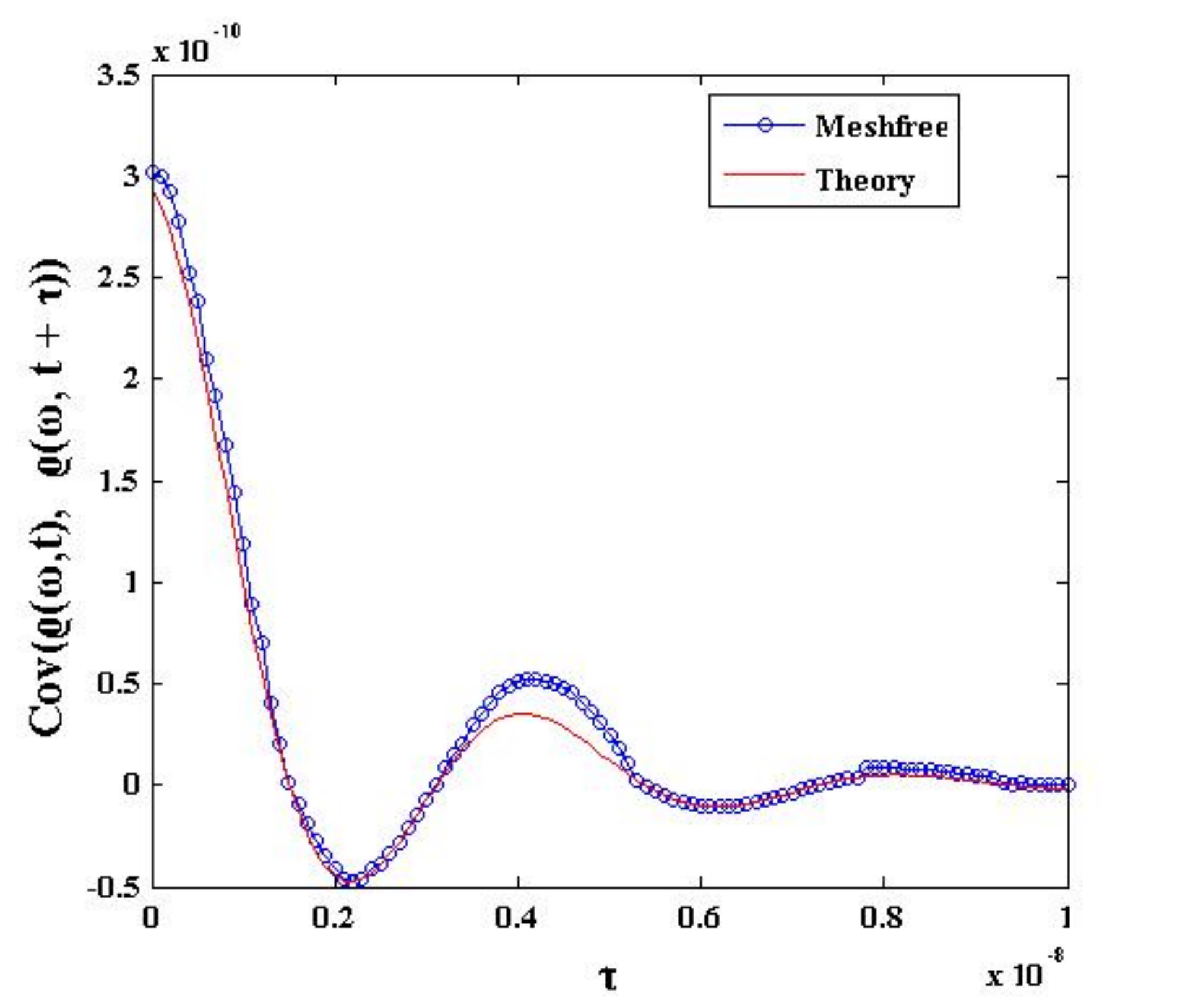}
\caption{Time covariance of density fluctuation for equilibrium problem on a periodic domain.}
\label{figure1}
\end{figure}

In figure \ref{figure1}, we compare the  theoretical covariance of the density with the described meshfree simulation. We find a reasonable agreement of the results up to time $  \left( \approx 4 \times 10^{-9}\right)$, when a sound wave crossed the system, i.e. the comparison with the theory is only accurate for short time because of the finite size effect.   
\subsection{Non-equilibrium}
In this last numerical test we consider a random walk of a standing shock wave due to spontaneous fluctuations. Our interest is the variance of the shock location as a function of time. The shock location is given by
\begin{equation}
\label{eq45}
\sigma_{\rho}\left( t \right) =  \int_{-L/2}^{\sigma(t)} \rho_{L}\, dx + \int^{L/2}_{\sigma(t)} \rho_{R}\, dx = \int_{-L/2}^{L/2} \rho(x, t)\, dx 
\end{equation}
\begin{equation}
\Longrightarrow \sigma_{\rho\left(t\right)} = L\dfrac{\bar{\rho}\left(t\right) - \frac{1}{2}\left(\rho_{L} + \rho_{R}\right)}{\rho_{L} - \rho_{R}}
\end{equation}
where, $\bar{\rho}$ is the instantaneous average density. 

The system parameters for the simulation of the standing shock is given below in the Table \ref{table4}. 

Mass density, temperature and velocity on both side of shocks are given by the Rankine-Hugoniot conditions. The boundary condition is adapted to the  initial condition. 
\begin{table}[h]
\begin{tabular}{ | p{3 cm} | l |l| p{3 cm} | l | }
\hline 
System length & $5.0 \times 10^{-4}$ && Reference Mean free path & $6.26 \times 10^{-6}$ \\ [0.2 cm]
System volume & $7.84 \times 10^{-16}$ && Time step & $1.0 \times 10^{-13}$ \\  [0.2 cm] \hline
RHS mass density & $1.78 \times 10^{-3}$ && LHS mass density & $4.07 \times 10^{-3}$ \\ [0.2 cm]
RHS velocity & $- 61562$ && LHS velocity& $ - 26933 $ \\ [0.2 cm]
RHS sound speed & $30781$ && LHS sound speed & $ 44373 $ \\ [0.2 cm]
RHS temperature & $273.0$ && LHS temperature & $567.0$ \\[0.2 cm] 
\hline
\end{tabular}
\caption{System parameter in CGS units for simulation of standing shock waves, Mach number 2.0}
\label{table4}
\end{table}

We considered two different shock strengths by considering  $\textrm{Mach 2.0 , Mach 1.4}$. We have compared our results with the  Finite Volume-third order Runge-Kutta (FVRK3) scheme in \cite{Garcia2}, which already has good agreement with DSMC molecular simulation for prescribed Mach numbers. 
\newpage
\begin{figure}[h]
\label{figure2}
\centering
\includegraphics[width=8cm,height=7cm]{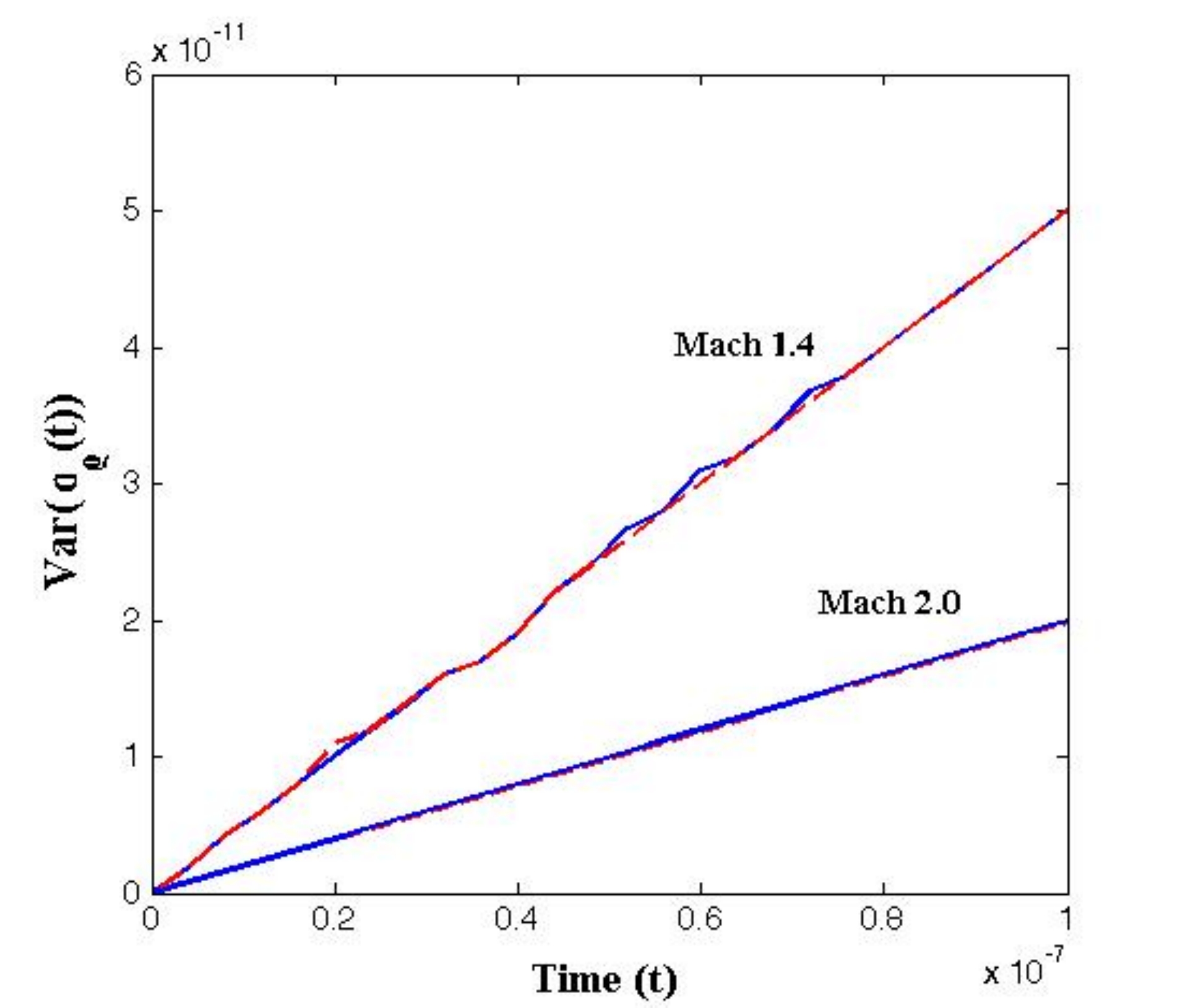}
\caption{variance of shock location. Solid lines represent Meshfree simulation and dashed lines are FVRK3. }
\end{figure}
\section{Concluding Remarks}

The results of the simulation show that the Lagrangian particle scheme gives a good agreement with the  theoretical values of variances for the conserved variables. This  shows that the scheme is able to accurately represent fluctuations in equilibrium flow. Further tests for the standing shock waves confirm that with a meshfree discretization we were able to reproduce stochastic drift of shock waves, as verified by comparison with the FVRK3 scheme from \cite{Garcia2}, which has been compared with molecular simulation. 

It has already been mentioned in earlier literature that the ability of continuum model to accurately capture the fluctuation is very much sensitive to the construction of the numerical scheme. This is also true for the meshfree framework.  Minor changes in implementation lead to significant changes in accuracy and behavior. In  future work the above will be extended to higher dimension. We will further study higher dimensional incompressible flow and the dynamics of small particles at fluid interfaces.

\bibliographystyle{model1-num-names}
\begin{small}
 
 \end{small}
\end{document}